\documentclass[a4paper,12pt]{article}
\def\beq{\begin{equation}}
\def\eeq{\end{equation}}
\def\bea{\begin{eqnarray}}
\def\eea{\end{eqnarray}}
\def\nn{\nonumber}
\def\bra#1{\left\langle #1\right|}
\def\ket#1{\left| #1\right\rangle}
\def\braket#1#2{\left\langle #1 \right| \left. #2 \right\rangle}
\def\sl2{{\cal U}_h(sl(2))}
\def\hf{\frac{1}{2}}
\def\F{{\cal F}}
\newtheorem{lemma}{Lemma}[section]
\newtheorem{prop}{Proposition}[section]
\newtheorem{thm}{Theorem}[section]
\newtheorem{DEF}{Definition}[section]
\newtheorem{cor}{Corollary}[section]
\usepackage{latexsym}

\begin{document}
\thispagestyle{empty}
\begin{flushright}
OWUAM-030 \\
March 17, 1999
\end{flushright}
\begin{center}

\vspace*{3cm}
{\LARGE \bf Symplecton for $ \sl2 $ \\ and Representations of $ SL_h(2) $}

\vspace{2.5cm}

{\large N. Aizawa

\bigskip
\textit{Department of Applied Mathematics \\ 
Osaka Women's University \\
Sakai, Osaka 590-0035, JAPAN}
}
\end{center}

\vfill
\begin{abstract}
  Polynomials of boson creation and annihilation operators which form irreducible tensor operators for Jordanian quantum 
algebra $ \sl2 $, called $h$-symplecton, are 
introduced and  
their properties are 
investigated. It is shown that many properties of symplecton for 
Lie algebra $ sl(2) $ are extended to $h$-symplecton. 
The $h$-symplecton 
is also a basis of irreducible representation of $ SL_h(2) $ dual to $ \sl2. $
As an application of the procedure used to construct $h$-symplecton, 
we construct 
the representation bases of $ SL_h(2) $  
on the quantum $h$-plane. 
\end{abstract}

\newpage
\setcounter{page}{1}
%
%
%
%
\section{Introduction}

  It is no doubt that well-developed representation theories 
are necessary when we apply algebraic objects to physics. The 
simplest examples in quantum physics are angular momentum algebra 
$ su(2) $ and rotation matrices in 3 dimensional space $ SO(3). $ 
Their algebraic structure is simple but contents of representation 
theories are quite rich \cite{BieL8}. To investigate these algebraic 
objects or their complexification could be a foundation for further 
investigation of higher dimensional objects. 

  As for deformation of Lie groups and Lie algebras, $q$-deformation 
of Lie algebra $ sl(2) $ and Lie group $ SL(2) $ (and their real form) 
is studied quite well. Their representation theories have attracted much interest 
in both physics and mathematics and give a way to higher dimensional cases 
\cite{BLohe1}. There exists, however, some other deformation of Lie 
groups and algebras and these are generally called nonstandard deformation.
The most studied one may be the so-called Jordanian deformation obtained by 
Drinfeld twist from a Lie algebra or a known quantum algebra. 
The simplest examples are, of course, the Jordanian deformation of Lie algebra 
$ sl(2) $ and its dual. The Jordanian deformation of Lie group $ SL(2) $,
denoted by $ SL_h(2) $, is 
studied in \cite{dem,zak,ewen} and then 
Ohn introduced its dual algebra, namely, Jordanian deformation of 
$ sl(2) $ denoted by $ \sl2 $ \cite{ohn}. The Jordanian quantum algebra 
$ \sl2 $ is more natural than the $q$-deformed $ sl(2) $ in the sense that 
it is regarded as the angular momentum algebra with nonstandard coproduct 
(see \S 3) and we can use ordinary boson operators to represent $ \sl2 $, 
while it is hard to regard the $q$-deformed $ sl(2) $ as angular momentum 
and $q$-deformed boson algebras are used for representations
\footnote[2]{There exist mappings from the ordinary boson algebra to 
$q$-deformed ones \cite{btoq,fiore,fiore2}. It is, however, simpler to use 
the $q$-deformed boson algebras for representation theories.}. 
However, the representation theories of $ \sl2 $ and $ SL_h(2) $ have 
not been developed yet. We do not know, for example, the Racha coefficients and 
matrix elements of the universal $R$-matrix for $ \sl2 $. As for $ SL_h(2) $, 
even its representation matrices are not obtained. 

  In this article, in order to develop representation theories for 
Jordanian deformed algebras, we study symplecton 
for $ \sl2 $ and apply it to investigate representation matrices of $ SL_h(2) $.
The use of symplecton could be legitimated by recalling the properties of 
symplecton and $q$-deformed case.
The symplecton, introduced by Biedenharn and Louck \cite{BieL,BieL2}, 
is a polynomial of boson creation and annihilation operators 
which form an irreducible 
tensor operator of $ sl(2) $, that is, symplecton is a basis of irreducible 
representation (irrep.) for both $ sl(2) $ and $ SL(2). $  It is known that the 
symplecton is written in terms of Gauss hypergeometric function and product of 
two symplecton is reduced to a series of symplecton with Racha coefficients. 
In Ref.\cite{BieL}, application of symplecton to the Elliot model for 
nuclei is discussed, then it is found that Weyl-ordered polynomials for 
position and momentum operators are equivalent to symplecton \cite{LBL}. 
Many properties of symplecton are inherited from $ sl(2) $ to $q$-deformed 
case \cite{BLohe1,BLbook,Nom}. The $q$-deformed symplecton, 
called $q$-symplecton, is a 
irreducible tensor operator so that it is a irrep. basis for 
$q$-deformed $sl(2)$ and $ SL(2) $. The $q$-symplecton is written in 
terms of $q$-hypergeometric function and product of two $q$-symplecton 
is reduced to a series of $q$-symplecton with $q$-Racha coefficients.  $q$-Deformation of Weyl-ordered polynomial 
\cite{GelF} is formulated with $q$-symplecton. 
These facts show that symplecton is a powerful tool to investigate 
representation. 

  The plan of this article is as follows. Next three sections are mainly 
preparation for symplecton of $ \sl2 $. We often call the symplecton for $ \sl2 $ 
$h$-symplecton. The next section is a review of symplecton for $ sl(2). $ 
Some of the properties of symplecton listed in \S 2 will be extended to 
$h$-symplecton. \S 3 is devoted to the Jordanian quantum algebra $ \sl2 $ and 
Jordanian quantum group $ SL_h(2) $. We give new results on the twist element 
and Racha coefficients for $ \sl2. $ In \S 4, tensor operators for a Hopf 
algebra is introduced according to Ref.\cite{rs} and the relation between 
tensor operators for a Lie algebra and a Hopf algebra obtained by 
Drinfeld twist is discussed. Applying the result in \S 4, the $h$-symplecton is 
constructed from the $ sl(2) $ symplecton in \S 5. The properties of 
$h$-symplecton are studied in \S 5 and \S 6. We shall consider 
another irreducible tensor operators 
obtained from the quantum $h$-plane for $ \sl2 $ in \S 7 and using 
these tensor operators, as well as $h$-symplecton, irreps. of 
$ SL_h(2) $ are considered. \S 8 is  concluding remarks.

%
%
%
%
\section{Symplecton for $ sl(2) $}

  The symplecton realization of $ sl(2) $ is said to be "minimal", 
since only one kind of boson operator is used. It is in marked contrast 
to the well-known Jordan-Schwinger realization where two kinds of bosons 
are necessary. Let us first review the definition and important properties 
of the $ sl(2)$ symplecton \cite{BieL, BieL2}.

  Let $ \bar a, a $ be boson operators satisfying $ [\bar a,\; a] = 1$, 
and define
\beq
 J_+ = -\frac{1}{2} a^2, \qquad J_- = \frac{1}{2} \bar a^2, 
  \qquad
  J_0 = \frac{1}{2}(a \bar a + \bar a a). 
  \label{sl2real}
\eeq
It is easy to verify that (\ref{sl2real}) satisfies the $ sl(2) $ commutation 
relations
\beq
   [J_0,\; J_{\pm}] = \pm 2 J_{\pm}, \qquad 
   [J_+,\; J_-] = J_0.               \label{sl2comm}
\eeq
The symplecton is a polynomial in $ \bar a$ and $ a $ and 
form a irreducible tensor operator of $ sl(2) $ belonging to the 
spin $j$ representation ($ j = \hf, 1, \frac{3}{2}, \cdots $). Namely 
the symplecton, denoted by $ P_j^m(a, \bar a), $ is defined by
\bea
 & & 
 [J_{\pm},\; P_j^m] = \sqrt{(j\mp m) (j\pm m +1)} P_j^{m\pm 1}, \nn \\
 & & 
 [J_0,\;   P_j^m] = 2m   P_j^m. \label{defsl2}
\eea

The basic idea of symplecton is to treat $ \bar a $ and $ a$ in a
symmetric way. To this end, the usual "boson calculus" is replaced with 
the so-called "symplecton calculus", that is, instead of the boson 
vacuum $ \ket{0} $ satisfying $ \bar a \ket{0} =0 $, 
the formal ket $ \ket{\ } $ 
which is not annihilated by both $ \bar a $ and $ a $ is introduced. 
The representation bases in the realization (\ref{sl2real}) are 
formed by letting $ P_j^m $ act on $ \ket{\ } $, and the action of 
generators on the bases is defined by 
$ J_{\alpha} \ket{jm} = [J_{\alpha},\; P_j^m] \ket{\ }$. 
There exists an appropriate definition of a inner product for these 
$ \ket{jm} $, so that we obtain the usual unitary representations of 
$ sl(2) $ with spin $j$. 

  The explicit form of the polynomials $ P_j^m(a, \bar a) $ is found by 
solving $ [J_+,\; P^j_j ] = 0 $ to obtain $ P_j^j = a^{2j} $, and then 
using the action of $ J_- $ to calculate $ P_j^m $, 
\beq
 P_j^m(a,\bar a) = 
              \frac{1}{2^{j-m}} \left[ \frac{(2j)! (j-m)!}{(j+m)!} \right]^{1/2} 
              \sum_{s=0}^{j-m}\;
              \frac{\bar a^{j-m-s} a^{j+m} \bar a^s}{s! (j-m-s)!}_.
    \label{exp1}
\eeq
An alternative form for $ P_j^m $ is obtained by starting with 
$ P_j^{-j} = \bar a^{2j} $ and 
then using the action of $ J_+ $,
\beq
  P_j^m(a,\bar a) =  
              \frac{1}{2^{j+m}} \left[ \frac{(2j)! (j+m)!}{(j-m)!} \right]^{1/2} 
              \sum_{s=0}^{j+m}\; 
              \frac{a^s \bar a^{j-m} a^{j+m-s}}{s! (j+m-s)!}_.
     \label{exp2}
\eeq

  We would like to list some properties of $ sl(2) $ symplecton. For their 
proof or detail, we refer the reader to Refs.\cite{BieL, BieL2}. 
  
\noindent
  (1) A set of polynomials $ \{P_j^m(a,\bar a) \ | \ m = -j, -j+1, \cdots, j\}$ 
forms representation bases for 
the Lie group $ SL(2) $ as well as the Lie algebra $ sl(2) $. The boson 
commutation relation is covariant under the action of $ SL(2) $ defined by
\beq
  (a', \; \bar a') = (a, \; \bar a) \left(
   \begin{array}{cc}
   x & u \\ v & y
   \end{array}
   \right)_,                               \label{actionSL2}
\eeq
where the $ 2 \times 2 $ matrix is an element of $ SL(2) $. 
The transformed polynomial $ P_j^m(a', \bar a') $ is decomposed into 
$ P_j^m(a, \bar a) $ multiplied by 
polynomials in the entries of $ SL(2) $ matrix.
\beq
 P_j^m(a', \bar a') = \sum_n  P_j^n(a,\bar a) d^j_{nm}(g).
  \quad
  g \in SL(2) \label{dfunc}
\eeq
The $ (2j+1) \times (2j+1) $ matrix $ d^j_{nm}(g) $ gives an irrep. of $ SL(2) $ 
and is called Wigner's $d$-function in terminology of physics.

\noindent
(2) The polynomials $ P_j^m(a, \bar a) $ have a generating function. Let 
$ \xi, \eta $ be ordinary c-numbers commuting with $ a, \bar a$. Then
\beq
 (\xi a + \eta \bar a)^{2j} = \sqrt{(2j)!} \sum_{m=-j}^j\; 
 \Phi_{jm}(\xi, \eta) P_j^m(a, \bar a),
 \label{genfun}
\eeq
where $ \Phi_{jm} $ are well-known representation bases of both $ sl(2) $ 
and $ SL(2) $,
\beq
 \Phi_{jm}(\xi,\eta) = \frac{\xi^{j+m} \eta^{j-m}}{\sqrt{(j+m)! (j-m)!}}_.
 \label{phi}
\eeq
Irreps. of $ sl(2) $ are constructed on (\ref{phi}) by the realization
\beq
  J_+ = \xi \frac{d}{d\eta}, \quad \ 
  J_- = \eta \frac{d}{d\xi}, \quad \ 
  J_0 = \xi \frac{d}{d\xi} - \eta \frac{d}{d\eta}, 
  \label{phisl2}
\eeq
while irreps. of $ SL(2) $ are obtained by the following transformation
\beq
   (\xi', \; \eta') = (\xi, \; \eta) \left(
   \begin{array}{cc}
   x & u \\ v & y
   \end{array}
   \right)_,                         \label{phiSL2}
\eeq
it follows that
\beq
    \Phi_j^m(\xi', \eta') =  \sum_n \Phi_j^n(\xi, \eta) \; d^j_{nm}(g),
    \label{phiSL2x}
\eeq
where we have obtained the same $d$-function as (\ref{dfunc}).

\noindent
(3) The symplecton polynomials can be expressed in terms of 
Gauss hypergeometric function $ {}_2F_1(a,b;c;z) $. The polynomial 
$ {}_2F_1(a,b;c;z) $ is defined by
\beq
 {}_2 F_1(a, b; c; z) = 
    \sum_{n=0}^{\infty}\; \frac{(a)_n (b)_n}{n! (c)_n}z^n, \label{2f1}
\eeq
where $ (a)_n $ stands for the sifted factorial
\beq
   (a)_n = \left\{
  \begin{array}{cl}
  1 & n = 0 \\
  a(a+1) \cdots (a+n-1) & n = 1, 2, \cdots
  \end{array}
  \right.                                      \label{siffac}
\eeq
Now define the operator $ N = a \bar a $, then the symplecton 
$ P_j^m(a, \bar a) $ is written in terms of $ {}_2F_1(a,b;c;z) $ with 
$ z = -1 $ and the parameters $ a,c $ become functions of operator $N$. 
The expression (\ref{exp1}) becomes
\bea
 P_j^m &=& \frac{1}{2^{j+m}}\left[
      \frac{(2j)!}{(j+m)!(j-m)!} \right]^{1/2} \frac{(N+j-m)!}{(N-2m)!} 
      \label{F21} \\
      & \times & {}_2F_1(-N+2m, -j+m; -N-j+m; -1) (\bar a)^{-2m}. \nn
\eea
In this way, properties of $ P_j^m $ are reduced to properties of the hypergeometric 
function. Especially, the equivalence of two form (\ref{exp1}) and (\ref{exp2}) 
is explained by the formula
\beq
  {}_2F_1(a, b; c; z) = (1-z)^{c-a-b} {}_2F_1(c-a, c-b; c; z).
  \label{F21formula}
\eeq

\noindent
(4) The polynomials $ P_j^m(a, \bar a) $ are transformed under the 
action $ a \rightarrow \bar a, \ \bar a \rightarrow -a $
\beq
  P_j^m(\bar a, -a) = (-1)^{j-m} P_j^{-m}(a, \bar a). \label{adj}
\eeq
To define an inner product for the bases $ \ket{jm} = P_j^m \ket{\ } $, 
the property (\ref{adj}) and the product formula discussed below play 
a crucial role.

\noindent
(5) Let $ P_j^m $ and $ P_{j'}^{m'} $ be the symplecton polynomials, then 
they obey the product law
\beq
   P_j^m  P_{j'}^{m'} = \sum_{k=|j-j'|}^{j+j'}\; 
  \bra{k}j\ket{j'} C_{m',m,m+m'}^{j',\; j,\quad \ k} 
  P_k^{m+m'},                                      \label{product}
\eeq
where 
\bea
  & & \bra{k}j\ket{j'} = 2^{k-j-j'} (2k+1)^{-1/2} \nabla(kjj'), \nn \\
  & & \nabla(abc) = \left[
      \frac{(a+b+c+1)!}{(a+b-c)! (a-b+c)! (-a+b+c)!} \right]^{1/2}_, \label{tri}
\eea
and $ C_{m',m,m+m'}^{j',\; j,\quad \ k} $ is the 
Clebsch-Gordan coefficient (CGC) for $ sl(2). $ 
The associativity of the products 
$ 
(P_a^{\alpha} P_b^{\beta}) P_c^{\gamma} 
 =
 P_a^{\alpha} (P_b^{\beta} P_c^{\gamma})
$
gives a relation between "triangle functions"
\beq
  \nabla(acf) \nabla(bdf) = (2f+1) \sum_e W(abcd;ef) \nabla(abe) \nabla(cde) ,
  \label{racha}
\eeq
where $ W(abcd;ef) $ is the Racha coefficient. 

  The inner product for $ \ket{jm} $ is defined by
\beq
 \braket{jm}{j'm'} = \bra{\ }(-1)^{j-m}P_j^{-m} \cdot P_{j'}^{m'} \ket{\ }, 
 \label{inner}
\eeq
and the operation $ \bra{\ }(\cdots) \ket{\ } $ means to take only the $ j= 0 $ 
part of the expression $ (\cdots) $. Applying the product law (\ref{product}) 
to the RHS of (\ref{inner}), we see that that the $j=0 $ part is given by 
the CGC $ C^{j,j',0}_{m,m',0} $, so that the bases $ \ket{jm} $ are 
orthonormal. 

%
%
%
%
\setcounter{equation}{0}
\section{Jordanian Deformation of $ sl(2) $ and $ SL(2) $}

 The Jordanian quantum algebras $ {\cal U}_h({\bf g}) $ are obtained 
from the (universal enveloping algebra ${\cal U}({\bf g}) $ of ) Lie algebras $ {\bf g} $ from
Drinfeld twist \cite{dri}. We denote the coproduct, conuit and 
antipode for $ {\cal U}({\bf g}) $, when it is regarded as a Hopf algebra, 
by $ \Delta, \epsilon, S $, respectively. 
With the invertible element 
$ {\cal F} \in {\cal U}({\bf g}) \otimes {\cal U}({\bf g}) $ satisfying
\bea
 & & (\epsilon \otimes id)(\F) = (id \otimes \epsilon)(\F) = 1, 
 \label{twist1} \\
 & & \F_{12} (\Delta \otimes id)(\F) = \F_{23}(id \otimes \Delta)(\F),
 \label{twist2}
\eea
the algebra $ {\cal U}_h({\bf g}) $ is defined by the same 
commutation relations as $ {\bf g} $ 
and the following Hopf algebra mappings
\beq
 \tilde \Delta = \F \Delta \F^{-1}, \qquad 
 \tilde \epsilon = \epsilon, \qquad
 \tilde S = u S u^{-1},                           \label{defHopf}
\eeq
where $ u = m(id \otimes S)(\F), \ u^{-1} = m(S \otimes id)(\F^{-1}) $, 
$ m $ denotes the usual product in $ {\bf g}. $  This is a triangular 
Hopf algebra whose universal $R$-matrix is given by $ {\cal R} = \F_{21} \F^{-1}. $

For the case of $ {\bf g} = sl(2) $, $ \F $ is given by \cite{ks}
\beq
 \F = \exp\left(-\hf J_0 \otimes \sigma \right), 
 \qquad
 \sigma = -\ln(1 -2hJ_+),                        \label{sigma}
\eeq
The twist element $ \F $ used here gives different form of $ \sl2 $ from 
the one in Ref.\cite{ohn}. The relationship between these two form 
is given in Appendix A. 
The explicit form of Hopf algebra mappings for $ \sl2 $ is summarized in  
Appendix B (some of them will be used in the later computation). 
An application of the $ \sl2 $ to the Heisenberg spin chain is found 
in Ref. \cite{ks}.
  The finite dimensional highest weight irreps. for $ \sl2 $ are same 
as $ sl(2) $, because of the same commutation relations.
We shall use the following lemmas on tensor product representations 
in subsequent sections.
\begin{lemma}{\rm \cite{ks}}
 Let $ V^{j_1}, V^{j_2} $ be the representation space with the highest 
weight $ j_1, j_2 $. Then the tensor product of them is completely reducible,  
$i.e.$
\[
 V^{j_1} \otimes V^{j_2} = \bigoplus_{j=|j_1-j_2|}^{j_1+j_2} V^j,
\]
and the bases of $ V^j $ are given by
\beq
 e^{(j_1j_2)j}_m = \sum C^{\; j_1,\; j_2,\; j}_{m_1,m_2,m} 
   F^{\ \ \ j_1,j_2}_{k_1,k_2,\ m_1,m_2} e^{j_1}_{k_1} \otimes e^{j_2}_{k_2},
   \label{tensor}
\eeq
where $  C^{\; j_1,\; j_2,\; j}_{m_1,m_2,m} $ is the 
CGC of $ sl(2) $ and 
$  F^{\ \ \ j_1,j_2}_{k_1,k_2,\ m_1,m_2} $ is the matrix element of 
$ \F $ on $ V^{j_1} \otimes V^{j_2}. $ 
\label{lem1}
\end{lemma}
The explicit form of matrix elements $  F^{\ \ \ j_1,j_2}_{k_1,k_2,\ m_1,m_2} $ 
is given in Appendix C (It seems to be the first time to show the 
explicit form of $  F^{\ \ \ j_1,j_2}_{k_1,k_2,\ m_1,m_2} $ in the literature,
and this also gives the explicit form of the $R$-matrix for $ \sl2 $).
\begin{lemma}
 The Racha coefficients for $ sl(2) $ and $ \sl2 $ coincide.
\label{lem2}
\end{lemma}
Lemma \ref{lem2} is proved in Appendix D.

  The matrix quantum group dual to $ \sl2 $ is called the Jordanian quantum 
group $ SL_h(2) $. It is generated by four elements $ x, y, u $ and $ v $ 
subject to the relations \cite{dem,zak,ewen}
\bea
 & & [v,\; x] = hv^2, \qquad [u,\; x] = h(1-x^2), 
     \nn \\
 & & [v,\; y] = h v^2, \qquad [u,\; y] = h(1-y^2),
     \label{SLh2} \\
 & & [x,\; y] = h(xv - yv), \qquad [v,\; u] = h(xv + vy).
     \nn
\eea
It follows that the central element of $ SL_h(2) $ 
which gives the determinant of the quantum 
matrix 
\beq
 T = \left(
  \begin{array}{cc}
  x & u \\ v & y
  \end{array}
  \right)_,      \label{quantumM}
\eeq
is defined by
\beq
detT = xy - uv - h xv = 1.     \label{det} 
\eeq
The $ SL_h(2) $ has a Hopf algebra structure. 
The relations (\ref{SLh2}) and Hopf 
algebra mappings are summarized in the FRT-formalism \cite{frt} 
with the $R$-matrix
\beq
  R = \left(
  \begin{array}{rrrr}
  1 & h & -h & h^2 \\
  0 & 1 & 0 &h \\
  0 & 0 & 1 & -h \\
  0 & 0 & 0 & 1
  \end{array}
  \right)_.                           \label{SLh2R}
\eeq
The coproduct, the counit and the antipode are given by
\bea
 & & \Delta(T) =  T \stackrel{\cdot}{\otimes} T, \nn \\
 & &  \epsilon(T) = \left(
  \begin{array}{cc}
  1 & 0 \\ 0 & 1
  \end{array}
  \right)_,              \nn \\
 & &  S(T) = T^{-1} = \left(
 \begin{array}{ccc}
 y - hv & \ \ &-u -h(y-x) + h^2 v \\
 -v & & x+hv
 \end{array}
 \right)_. \nn
\eea

  Let us define the $d$-function for $ SL_h(2) $ using the notion of 
comodule. A vector space $M$ is called right $SL_h(2)$ comodule 
if there is a map $ \rho : M \rightarrow M \otimes SL_h(2) $ such that 
the following relations are satisfied
\beq
  (\rho \otimes id)\circ \rho = (id_M \otimes \rho) \circ \rho, 
  \qquad
  (id_M \otimes \epsilon) \circ \rho = id_M,
  \label{comodule}
\eeq
where $ id_M $ stands for the identity map in $M$. Using bases $ e_i $ of 
$ M$, the map $ \rho $ is written as
\beq
  \rho(e_i) = \sum_j e_j \otimes \tilde d_{ji}, \label{dS1}
\eeq
it follows that the relations (\ref{comodule}) are rewritten as
\beq
  \Delta(\tilde d_{ij}) = \sum_k \tilde d_{ik} \otimes \tilde d_{kj}, 
  \qquad
  \epsilon(\tilde d_{ij}) = \delta_{ij}.
  \label{dS2}
\eeq
We call the $ \tilde d_{ij} $ satisfying (\ref{dS1}) and (\ref{dS2}) 
the $d$-function for $ SL_h(2) $. In the following sections, we 
deal with the case in which the vector space $M$ has an algebraic structure. 
It is natural, in this case, to require that the map $ \rho $ 
should respect the extra structure on $M$. 

%
%
%
%
\setcounter{equation}{0}
\section{Tnsor Operators and Twist}

  To define the symplecton for $ \sl2$, it is necessary to 
extend the notion of tensor operators to Hopf algebra. This has 
been carried out by Rittenberg and Scheunert \cite{rs}. 
Tensor operators are defined for each realization of the Hopf 
algebra $ {\cal H} $ under consideration. Assuming that we have a realization 
of $ {\cal H} $, we first define the adjoint action.
\begin{DEF}
 Let $ W, W'$ be representation space of $ {\cal H} $, and let $t$ be an 
operator which carries $ W $ into $W'$. Then the adjoint action of $ X \in {\cal H} $ 
on $t$ is defined by
\beq
    adX(t) = m(id \otimes S)(\Delta(X)(t \otimes 1)).  \label{adj2}
\eeq
\end{DEF}
The adjoint action has two important properties
\beq
  ad XX'(t) = ad X \circ ad X'(t), \qquad
  ad X(t \otimes s) = \sum_i ad X_i(t) \otimes ad X'_i(s), \label{adjpro}
\eeq
where the coproduct for $X$ is written as 
$ \Delta(X) = \sum_i X_i \otimes X'_i.$ 
From these properties, we see that the adjoint action gives a 
representation of $ {\cal H} $
\beq
  ad[X,\; X'] (t) = [ad X,\; ad X'](t).
\eeq
Tensor operators for $ {\cal H} $ are defined as operators which form 
representation bases 
of $ {\cal H} $ under the adjoint action. 
\begin{DEF}
  Let $ D(X) $ be a 
representation matrix of $ X \in {\cal H} $. The operators $ t_{\alpha} $ 
are called the tensor operator, if they satisfy the relation
\beq
  ad X(t_{\alpha}) = \sum_{\beta} D(X)_{\beta \alpha} t_{\beta}. \label{top}
\eeq
If the representation 
is irreducible, the tensor operators are called irreducible tensor operators.
\end{DEF}

 The explicit form of the adjoint action for $ \sl2 $ reads
\bea
   & & adJ_0(t) = [J_0, \; t] e^{-\sigma}, \nn \\
   & & adJ_+(t) = e^{-\sigma}[J_+e^{\sigma},\; t], \label{adjsl2} \\
   & & adJ_-(t) = [J_- + hJ_0 + \frac{h}{2}J_0^2,\; t]e^{-\sigma}
       -h[J_0,\; t]e^{-2\sigma}
       -\frac{h}{2}[J_0,\; [J_0,\; t]]e^{-2\sigma}. \nn 
\eea
Some examples of the $ \sl2 $ tensor operators are considered in Ref.\cite{na1} 
and they are applied to construct boson algebra which is covariant 
under the action of Jordanian matrix quantum groups \cite{ques1}. 

 Since the coproduct for Lie algebra  and Jordanian quantum algebra is related via twist element 
(\ref{defHopf}), tensor operators for these algebras are also 
related by twisting via $ \F $ \cite{fiore}.
\begin{lemma}
Let $ t_{\alpha} $ be tensor operators for Lie algebra $ {\bf g} $ and 
$ \tilde t_{\alpha} $ be corresponding ones for Jordanian quantum 
algebra $ {\cal U}_h({\bf g}). $ Then these tensor operators are 
related via the twist element $ \F $
\bea
 & & \tilde t_{\alpha} = m(id \otimes \tilde S) (\F (t_{\alpha} \otimes 1) \F^{-1}), 
     \label{tt1} \\
 & & t_{\alpha} = m(id \otimes S) (\F^{-1} (\tilde t_{\alpha} \otimes 1) \F).
     \label{tt2}
\eea
\label{lem3}
\end{lemma}
\textit{Proof} : The first relation (\ref{tt1}) is derived in Ref.\cite{fiore} 
(Proposition 3). The second one (\ref{tt2}) is its inverse. The expression 
used in Lemma \ref{lem3} is different form Ref.\cite{fiore}, it may be good to 
show the second relation as an example of the proof. 
It is proved by showing the substitution of  (\ref{tt2}) into (\ref{tt1}) gives 
the identity map. 

  Let us write the twist element and its inverse as
\[
  {\cal F} = \sum f^a \otimes f_a, \qquad 
  {\cal F}^{-1} = \sum g^a \otimes g_a,
\]
then
\[
  u = \sum f^a S(f_a), 
 \quad
 u^{-1} = \sum S(g^a) g_a,
\]
and the relation (\ref{tt1}) becomes
\beq
 \tilde t_{\alpha} 
 = \sum f^a t_{\alpha} g^b \tilde S(f_a g_b) 
 = \sum f^a t_{\alpha} g^b u S(f_a g_b) u^{-1} 
 = \sum f^a t_{\alpha} S(f_a) u^{-1}, \label{one}
\eeq
where we used
\[
  \sum g^b u S(g_b) = \sum g^b f^a S(g_b f_a) = m(id \otimes S)(\F^{-1} \F) = 1.
\]
On the other hand, the relation (\ref{tt2}) is rewritten
\beq
 t_{\alpha} = \sum g^a \tilde t_{\alpha}  f^b S(g_a f_b) 
 = \sum g^a \tilde t_{\alpha} u S(g_a). 
 \label{two}
\eeq
Substituting (\ref{two}) into (\ref{one})
\[
 \tilde t_{\alpha} = \sum f^a g^b \tilde t_{\alpha} u S(f_a g_b) u^{-1}
  = \sum f^a g^b \tilde t_{\alpha} \tilde S(f_a g_b) 
  = m(id \otimes \tilde S)(\F \F^{-1} (\tilde t_{\alpha} \otimes 1)) 
  = \tilde t_{\alpha}.
\]
This proves the second relation in Lemma \ref{lem3}. \hfill $\Box$

%
%
%
%
\setcounter{equation}{0}
\section{Symplecton Polynomials for $ \sl2 $}

  In this section, we derive the explicit form of the symplecton for 
$ \sl2 $ and investigate its properties. Since $ \sl2 $ has the 
same commutation relations as $ sl(2) $, $ \sl2 $ and $ sl(2) $ 
have the same realizations. 
Therefore the symplecton realization for $ \sl2 $, 
which is identical to the one for $ sl(2),$  
is the realization in terms of the usual boson operators. 
This is a contrast to the $q$-symplecton where the $q$-deformed boson 
operators are used. 

 Let $ \bar a $ and $ a $ be boson operators satisfying $ [\bar a, \; a] = 1 $, 
then the generators of $ \sl2 $ are realized by
\beq
  J_+ = -\frac{1}{2} a^2, \quad J_- = \frac{1}{2} \bar a^2, 
  \quad
  J_0 = \frac{1}{2}(a \bar a + \bar a a),          \label{slh2real}
\eeq
The $h$-symplecton, denoted by $ \tilde P_j^m(a, \bar a) $, is 
defined as a polynomial in $ \bar a, a $ satisfying
\bea
 & & ad J_{\pm}(\tilde P_j^m) = \sqrt{(j\mp m)(j\pm m+1)} \tilde P_j^{m\pm 1}, 
 \nn \\
 & & ad J_0(\tilde P_j^m) = 2m \tilde P_j^m, \label{defhsym}
\eea
where the adjont action on the LHS is given by (\ref{adjsl2}).
Using Lemma \ref{lem3}, 
the explicit form of $h$-symplecton is obtained from the corresponding one 
for $ sl(2). $
\begin{prop}
The explicit form of the $h$-symplecton defined by (\ref{defhsym}) is given by
\beq
  \tilde P_j^m(a,\bar a) = P_j^m(a,\bar a) e^{m \sigma},  \label{hsymp}
\eeq
where $ \sigma $ is given in (\ref{sigma}) and $ P_j^m(a, \bar a) $ 
denotes $ sl(2) $ symplecton. 
\label{prop1}
\end{prop}
\textit{Proof} : By definition of $ sl(2) $ symplecton, it holds that
\[
  (J_0 - 2m) P_j^m = P_j^m J_0.
\]
Using this and the RHS of (\ref{one}), 
\[
 \tilde P_j^m = \sum_{n=0}^{\infty} \frac{1}{n!} \left(-\frac{1}{2}\right)^n 
 P_j^m (J_0 + 2m)^n S(\sigma)^n u^{-1}
 = P_j^m
 \sum_{n=0}^{\infty} \sum_{s=0}^n \frac{(-1)^n (2m)^s}{2^n (n-s)! s!} 
 J_0^{n-s} S(\sigma)^n u^{-1}.
\]
Changing the order of sum, then replacing $n-s$ with $n$, we obtain
\beq
 \tilde P_j^m = P_j^m \sum_{s,n = 0}^{\infty} 
 \frac{(-1)^{n+s} (2m)^s}{2^{n+s} n! s!} J_0^n S(\sigma)^{n+s} u^{-1}.
 \label{three}
\eeq
Note that
\[
  u = \sum_{n=0}^{\infty} \left(-\frac{1}{2}\right)^n \frac{1}{n!} J_0^n S(\sigma)^n,
\]
and (\ref{defHopf}), it follows that (\ref{three}) is rewritten as
\[
 \tilde P_j^m = P_j^m \sum_{s=0}^{\infty} \left(-\frac{1}{2}\right)^s 
 \frac{(2m)^s}{s!} \tilde S(\sigma)^s 
 = P_j^m e^{m\sigma},
\]
where (\ref{sigmaHopf}) is used in the last equality. 
\hfill $\Box$

  We would like to show some explicit form of $h$-symplecton. 
For $ j = 1/2 $
\beq
 \tilde P_{1/2}^{-1/2} = \bar a e^{-\sigma/2} \equiv \bar a_h,\quad
 \tilde P_{1/2}^{1/2} = \bar a e^{\sigma/2} \equiv a_h,
 \label{jhalf}
\eeq
and for $ j = 1 $
\bea
 & & \tilde P_1^{-1} = \bar a^2 e^{-\sigma} = \bar a_h^2 + h\bar a_h a_h, \nn \\
 & & \tilde P_1^{0} = (\bar a a + a \bar a)/\sqrt{2} 
     = (\bar a_h a_h + a_h \bar a_h - ha_h^2)/\sqrt{2}, \label{jone} \\
 & & \tilde P_1^{1} = a^2 e^{\sigma} = a_h^2. \nn
\eea
The $ j=1/2 $ $h$-symplecton forms covariant $h$-deformed oscillator algebra
\beq
   [\bar a_h, \; a_h ] = 1 - h a_h^2,      \label{hosc}
\eeq
$i.e.$, the commutation relation (\ref{hosc}) is preserved under the action 
of $ SL_h(2) $
\beq
  (a_h', \; \bar a_h') = (a_h, \; \bar a_h) \left(
   \begin{array}{cc}
   x & u \\ v & y
   \end{array}
   \right)_.                        \label{action}
\eeq
This shows that it is possible to construct representations of $ SL_h(2) $ 
on $h$-symplecton. We shall discuss it later. 
It may be worth noting that the action (\ref{action}) is different from the 
ones in \cite{ques1,fiore} where $ a $ and $ \bar a $ are \textit{not} 
mixed by the action of quantum groups.

 The $ j=1 $ $h$-symplecton forms 
an algebra isomorphic to $ sl(2) $. Its commutation relations are
\bea
  & & [ P_1^0,\; P_1^1] = 2\sqrt{2} P_1^1 (1 - hP_1^1), \nn \\
  & & [ P_1^0,\; P_1^{-1}] = -2\sqrt{2} P_1^{-1}(1 - h P_1^1), \label{j1symp} \\
  & & [ P_1^1\; P_1^{-1}] = -2\sqrt{2} (1 - h P_1^1) P_1^0. \nn
\eea
The generators of $ sl(2) $ are written in terms of $ P_1^m $
\beq
   J_+ = -\frac{1}{2} P_1^1 (1-hP_1^1), \qquad
   J_0 = \frac{1}{\sqrt 2} P_1^0, \qquad
   J_- = \frac{1}{2} P_1^{-1} (1 - h P_1^1).    \label{JbyP}
\eeq

  We see, from the explicit form of $h$-symplecton (\ref{hsymp}), that 
the $h$ dependence of polynomial $ \tilde P_j^m(a, \bar a) $ is absorbed 
in $\sigma $ which is a infinite polynomial in $ a^2 $. 
Recall that the relationship 
between $ sl(2) $ symplecton and Gauss hypergeometric function $ {}_2F_1 $ 
is given in terms of the operator $ N = a \bar a $, then we see that 
the factor $ e^{m\sigma} $ in (\ref{hsymp}) does not affect 
this relationship. Therefore the specific hypergeometric 
function for $h$-symplecton may be again $ {}_2F_1 $. 

  The fact that the $j=1/2$ $h$-symplecton forms covariant $h$-oscillator algebra 
may suggest that it is useful to write $h$-symplecton in terms of covariant 
$h$-oscillators (\ref{jhalf}). 
\begin{prop}
The $h$-symplecton is written in terms of covariant $h$-oscillators as 
follows.
The corresponding expression for (\ref{exp1}) is 
\bea
   \tilde P_j^m(a_h, \bar a_h) &=& 
   \frac{1}{2^{j-m}} \left[ \frac{(2j)! (j-m)!}{(j+m)!} \right]^{1/2} 
              \sum_{s=0}^{j-m}\; 
              \frac{1}{s! (j-m-s)!} \nn \\
   &\times & 
   \bar a_h ( \bar a_h + h a_h) \cdots \{ \bar a_h + (j-m-s-1) h a_h \} 
   a_h^{j+m} \label{expo1} \\
   &\times &
   \{ \bar a_h - (2m+s) h a_h \} \{ \bar a_h - (2m+s-1)ha_h \} 
   \cdots 
   \{\bar a_h - (2m+1) h a_h \},\nn
\eea
and for (\ref{exp2}) is
\bea
   \tilde P_j^m(a_h, \bar a_h)  &=& 
   \frac{1}{2^{j+m}} \left[ \frac{(2j)! (j+m)!}{(j-m)!} \right]^{1/2} 
              \sum_{s=0}^{j+m}\; 
              \frac{1}{s! (j+m-s)!} \label{expo2} \\
        &\times &      
        a_h^s (\bar a_h-hsa_h) \{\bar a_h + h(1-s)a_h\} \cdots 
        \{\bar a_h + h(j-m-1-s)a_h\} a_h^{j+m-s}. \nn
\eea
\label{prop2}
\end{prop}
\textit{Proof} : From (\ref{jhalf})
\[
  \bar a = \bar a_h e^{\sigma/2}, \qquad
  a = a_h e^{-\sigma/2}.
\]
Substituting these into (\ref{exp1}) and (\ref{exp2}), then straightforward 
calculation proves the proposition. \hfill $ \Box $

  In order to discuss generating functions for $h$-symplecton, 
it is possible to apply  Lemma \ref{lem3} to the generating 
function (\ref{genfun}) for $ sl(2) $ symplecton , 
since the RHS of (\ref{genfun}) is a sum of tensor operators of $ sl(2). $ 
It follows that the RHS of (\ref{genfun}) becomes the sum of $h$-symplecton ; 
$ \displaystyle{ \sqrt{(2j)!} \sum_{m=-j}^j \Phi_{jm} \tilde P_j^m }. $ 
However the LHS may be quite complicated and may not be in closed form. 
Another way to obtain generating functions for $h$-symplecton is to 
substitute (\ref{hsymp}) and (\ref{jhalf}) into (\ref{genfun}) 
\beq
 ( \xi a_h e^{-\sigma/2} + \eta \bar a_h e^{\sigma/2} )^{2j} 
 = 
 \sqrt{(2j)!} \sum_{m=-j}^j\; \Phi_{jm}(\xi, \eta) \tilde P_j^m(a_h, \bar a_h) e^{-m\sigma}.
 \label{hgenfun}
\eeq
It is possible to remove $ \sigma $ from (\ref{hgenfun}) by usint 
the relation
$
 e^{\sigma} + h a_h^2 =1 
$, however the obtained 
relation is quite complicated. Therefore the simplest generating function for 
$h$-symplecton may be (\ref{hgenfun}) where the $ \sigma $ is regarded as 
a independent quantity subject to the relations
\beq
  [\sigma,\; a_h] = 0, \qquad
  [\sigma,\; \bar a_h] = 2h a_h,   \label{sigah}
\eeq
and $ \displaystyle{\lim_{h \rightarrow 0}} \sigma = 0. $

%
%
%
%
\setcounter{equation}{0}
\section{Product Law for $h$-Symplecton}

  It is possible to extend the product law (\ref{product}) for 
$ sl(2) $ symplecton to $h$-symplecton. The product law plays a 
crucial role when the symplecton calculus is considered. In this 
section, we first prove the product law for $h$-symplecton by 
using the one for $ sl(2)$ symplecton, then consider the 
symplecton calculus for $ \sl2. $
\begin{thm}
Let $ \tilde P_j^m $ and $ \tilde P_{j'}^{m'} $ be $h$-symplecton, 
then these obey the product law
\beq
  \tilde P_j^m \tilde P_{j'}^{m'} = 
 \sum_{k=|j-j'|}^{j+j'} \sum_{n,n'}\; \bra{k} j \ket{j'} 
 (F^{-1})_{n,n' \ m,m'}^{\quad j,j'} 
 C_{n', m, n'+m}^{j',\; j,\ \ \ k} \tilde P_{k}^{n'+m},
 \label{hproduct}
\eeq
where $ C_{m_1,m_2,m}^{j_1,j_2,j} $ is CGC for $ sl(2) $ and 
\bea
  & & \bra{k}j\ket{j'} = 2^{k-j-j'} (2k+1)^{-1/2} \nabla(kjj'), \nn \\
  & & \nabla(abc) = \left[
      \frac{(a+b+c+1)!}{(a+b-c)! (a-b+c)! (-a+b+c)!} \right]^{1/2}_. \label{tri2}
\eea
\label{thm1}
\end{thm}
\textit{Proof} : From Proposition \ref{prop1},
\[
 \tilde P_j^m \tilde P_{j'}^{m'} = P_j^m e^{m\sigma} \tilde P_{j'}^{m'}
 e^{-m\sigma} e^{m\sigma} . 
\]
Using the Hopf algebra mappings for $ \sigma $ given in (\ref{sigmaHopf}), 
we see that the adjoint action of $ e^{m \sigma} $ is given by
\beq
  adj e^{m \sigma}(t) = e^{m \sigma} t e^{-m \sigma}.
  \label{adjsigma}
\eeq
$h$-Symplecton is a irreducible tensor operator of $ \sl2 $, it follows that
\bea
  e^{m \sigma} \tilde P_{j'}^{m'} e^{-m \sigma} &=& 
  adj e^{m \sigma}(\tilde P_{j'}^{m'}) 
  = \sum_{n'}\; (e^{m \sigma})^{\; j'}_{n',m'} \tilde P_{j'}^{n'} \nn \\
  &=& \sum_{n,n'}\; \delta_{n,m} (e^{m \sigma})^{\; j'}_{n',m'} \tilde P_{j'}^{n'} \nn \\
  &=& \sum_{n,n'}\; (F^{-1})_{n,n' \ m,m'}^{\quad j,j'} \tilde P_{j'}^{n'}, 
  \label{stilPs}
\eea
where the matrix elements of $ \F $ (\ref{rhs}) is used in the last equality. 
Therefore, we have
\[
  \tilde P_j^m \tilde P_{j'}^{m'} = 
  \sum_{n,n'}\; (F^{-1})_{n,n' \ m,m'}^{\quad j,j'} P_j^m \tilde P_{j'}^{n'}
  e^{m \sigma}
  = \sum_{n,n'}\; (F^{-1})_{n,n' \ m,m'}^{\quad j,j'} P_j^m P_{j'}^{n'}
  e^{(n'+m) \sigma}
\]
Applying the product law (\ref{product}) for $ sl(2) $ symplecton, 
Theorem \ref{thm1} is proved.
\hfill $\Box $

\begin{cor}
The associativity of the products 
$
 (\tilde P_a^{\alpha} \tilde P_b^{\beta}) \tilde P_c^{\gamma} 
 = \tilde P_a^{\alpha} (\tilde P_b^{\beta} \tilde P_c^{\gamma}) 
$
gives the same relation as (\ref{racha}) for the triangle function $ \nabla(abc) $ 
appeared in Theorem \ref{thm1}.
\label{cor1}
\end{cor}
\textit{Proof} : The associativity gives the same relation as (\ref{racha}),  
but the Racha coefficients are replaced with the ones for $ \sl2 $. 
From Lemma \ref{lem2}, these two kinds of Racha coefficients coincide.
\hfill $\Box $

  Let us now consider the $h$-symplecton calculus. We assume the formal 
ket $ \ket{\ } $ and that both $ \bar a \ket{\ } $ and $ a \ket{\ } $ are 
nonvanishing vectors. Then the vectors defined by 
$ \ket{jm} = \tilde P_j^m \ket{\ } $ are irrep. bases of $ \sl2 $ 
provided that the action of $ X \in \sl2 $ is defined by 
$ X \ket{jm} = adj X(\tilde P_j^m) \ket{\ }. $  The dual bases are 
defined by $ \bra{jm} = \bra{\ } \tilde P_j^{-m}(-1)^{j-m} $ in order to 
keep the correspondence with the $h=0 $ case. The action of $ X \in \sl2 $ is, 
of course, given by $ \bra{jm} = \bra{\ } adjX(\tilde P_j^{-m}) (-1)^{j-m} $. 
The inner product is defined in the same manner as $ h=0 $ case, namely, 
\[
  \braket{jm}{j'm'} = \bra{\ }(-1)^{j-m}\tilde P_j^{-m} \cdot \tilde P_{j'}^{m'} \ket{\ }, 
\]
the operation $ \bra{\ } ( \cdots ) \ket{\ } $ means to take only the $ j= 0 $ 
part of the expression $ (\cdots) $. Applying the product law (\ref{hproduct}) for 
$h$-symplecton, we obtain
\beq
   \braket{jm}{j'm'} = \delta_{j,j'} 2^{-2j} F_{-m,m\ -m,m'}^{\qquad j,j}.
   \label{hinner}
\eeq
Therefore the vectors $ \ket{jm} $ and $ \ket{j'm'} $ are orthonormal 
if they belong to different irreps. but $not$ orthonormal if they 
belong to a same irrep. The nonvanising part on the RHS of (\ref{hinner}) 
depends on only the twist element $ \F$. 

 From the product law, we can show the following relations for 
$h$-symplecton
\begin{prop}
The following relations hold for $h$-symplecton.
\bea
 & & \sum_{m,m'}\; \tilde P_j^m \tilde P_{j'}^{m'} 
     F^{\quad j,j'}_{m,m'\ \ell,\ell'} = 
     \sum_k \bra{k} j \ket{j'} C^{j',j,\ \ k}_{\ell',\ell,\ell'+k} \; 
     \tilde P_k^{\ell + \ell'},
     \label{hrel1} \\
 & & \tilde P_{j'}^{m'}(a_h, \bar a_h + 2hma_h) = 
     \sum_{\ell=0} (F^{-1})_{n,n'\ m,m'}^{\quad j,j'}\; 
     \tilde P_{j'}^{m'+\ell}(a_h, \bar a_h).
     \label{hrel2}
\eea
\label{prop3}
\end{prop}
\textit{Proof} : The relation (\ref{hrel1}) is easily proved by 
multiplying the product law (\ref{hproduct}) by 
$ F^{\quad j,j'}_{n,n'\ m,m'} $ and summing over $ m, m'$. 
The relation (\ref{hrel2}) is derived by 
moving $ e^{m\sigma} $ to the right of $ P_{j'}^{m'} $ in (\ref{stilPs}). 
One can do that by using the relations
\[
  e^{m\sigma} a_h = a_h e^{m\sigma}, \qquad
  e^{m\sigma} \bar a_h = (\bar a_h + 2hm a_h) e^{m\sigma}.
\]
\hfill $\Box$

%
%
%
%
\setcounter{equation}{0}
\section{Quantum $h$-Plane and Representations of $ SL_h(2) $}

  The Jordanian quantum algebra $ \sl2 $ and Jordanian quantum 
group $ SL_h(2) $ are dual each other. It follows that 
any representation basis of $ \sl2 $ is also representation basis 
for $ SL_h(2) $ belonging to the same representation. Since 
$h$-symplecton is a irrep. basis of $ \sl2 $, it is also 
a irrep. basis of $ SL_h(2) $. We have seen this for $ j=1/2 $ 
in \S 5. The relation (\ref{action}) can be generalized to 
arbitrary $j$
\beq
  \tilde P_j^m(a_h', \bar a_h') = 
  \sum_n \tilde P_j^n(a_h,\bar a_h)  \tilde d^j_{nm}(g)
  \quad
  g \in SL_h(2).   \label{dands}
\eeq
We can obtain $d$-functions for $ \sl2 $ by substituting (\ref{action}) 
into the explicit form of $h$-symplecton given in Proposition \ref{prop2}. 
However, as is seen from the explicit form, 
the actual computation seems to be complicated.

  The use of quantum $h$-plane \cite{kar} provides us 
a procedure which is a little bit simpler in computation. 
In this section, we shall find irrep. bases for $ SL_h(2) $ in terms of quantum 
$h$-plane which give the same irreps. as $h$-symplecton 
by using the tensor operator approach. 
Recall that the functions $ \Phi_{jm}(\xi, \eta) $ defined by 
(\ref{phi}) are irrep. bases of $ sl(2) $ in the realization 
(\ref{phisl2}) and irrep. bases of $ SL(2) $ under (\ref{phiSL2}) 
as well. We can regard $ \Phi_{jm}(\xi, \eta) $ as a 
irreducible tensor operator of $ sl(2) $, since it is easy 
to verify that
\bea
 & & [J_{\pm}, \;  \Phi_{jm}] = \sqrt{(j\mp m)(j\pm m+1)} \Phi_{j\; m\pm1}, \nn \\
 & & [J_0,\;  \Phi_{jm}] = 2m \Phi_{jm}. \label{phiastensor}
\eea
From Lemma \ref{lem3}, it ie easy to find the corresponding 
irreducible tensor operators for $ \sl2 $.
\begin{prop}
 Let $ \xi, \eta $ be commutative numbers, then the 
followings are irreducible tensor operators for $ \sl2 $
\beq
  \tilde \Phi_{jm}(\xi, \eta) = \Phi_{jm}(\xi, \eta) e^{m\sigma}, \label{phih}
\eeq
where $ \sigma = -\ln(1-2h\xi \frac{d}{d\eta}) $.
\label{prop4}
\end{prop}
For $ j=1/2$, we have
\beq
 \tilde \Phi_{\frac{1}{2}\; \frac{1}{2}} = \xi e^{\sigma/2} 
 \equiv \xi_h, 
 \qquad
 \tilde \Phi_{\frac{1}{2}\; -\frac{1}{2}} = \eta e^{-\sigma/2} 
 \equiv \eta_h, 
 \label{hplane}
\eeq
and they satisfy the commutation relation
\beq
  [\xi_h, \; \eta_h ] = h \xi_h^2, \label{hplanecomm}
\eeq
this corresponds to the commutation relation of quantum 
$h$-plane in Ref.\cite{kar}. 
It is easily verified that the commutation relation (\ref{hplanecomm}) 
is preserved under the action of $ SL_h(2) $
\beq
  (\xi_h', \eta_h') = (\xi_h, \eta_h) \left(
 \begin{array}{cc}
  x & u \\ v & y
 \end{array}
 \right)_.  \label{hpaction}
\eeq
It is a easy exercise to write $ \tilde \Phi_{jm} $ in terms of 
$ \xi_h $ and $\eta_h$. Then $ \tilde \Phi_{jm}(\xi_h, \eta_h) $ 
forms irrep. bases of 
$ SL(2) $, that is,  the $d$-functions 
for $ SL_h(2) $ are obtained by substituting (\ref{hpaction}) into 
$ \tilde \Phi_{jm}(\xi_h, \eta_h) $. 
\begin{prop}
Irreps. of $ SL_h(2) $ on the quantum $h$-plane are 
obtained by 
\beq
  \tilde \Phi_{jm}(\xi_h',\eta_h') = \sum_k \; \tilde \Phi_{jk}(\xi_h, \eta_h) 
 \tilde d^{j}_{km},
 \label{SLh2onplane}
\eeq
where the irrep. bases are given by
\bea
    \tilde \Phi_{jm} &=& c_{jm} \; \xi_h^{j+m} (\eta_h - h(j+m)\xi_h) 
  (\eta_h - h(j+m-1)\xi_h) \nn  \cdots (\eta_h - h(2m+1)\xi_h), 
  \nn \\
  &=& c_{jm} \; \eta_h (\eta_h + h\xi_h) \cdots (\eta_h+(j-m-1)h\xi_h) \; 
      \xi_h^{j+m}, \label{hphix}
\eea
with 
\[
  c_{jm} = \frac{1}{\sqrt{(j+m)! (j-m)!}}.
\]
\label{prop5}
\end{prop}

  Since $ \tilde P_{jm} $ and $ \tilde \Phi_{jm} $ give 
the same irreps. of $ \sl2 $, they also give the same 
$d$-functions of $ SL_h(2) $. Indeed, the explicit computation shows  
that we obtain the same $d$-functions for $ j = 1/2 $ and $ j=1 $. 
The $ j=1/2 $ case gives the $ 2 \times 2 $ quantum matrix  $T$ 
(\ref{quantumM}) itself, while 
$ j=1 $ $d$-function reads
\beq
  d^1 = \left(
 \begin{array}{ccc}
  x^2 + h xv & \sqrt{2}(ux + h uv) & u^2+hu(x+y+hv) \\
  \sqrt{2}xv & 1 + 2 uv & \sqrt{2}(uy + huv) \\
  v^2 & \sqrt{2}yv & y^2 + hyv
 \end{array}
 \right)_.  \label{j1d}
\eeq
The $d$-functions for $ SL_h(2) $ are also discussed in Ref.\cite{CQ} where the 
authors assert that the $d$-functions can be obtained from the q-deformed ones 
via a contraction method and show some explicit examples. 
Another way to obtain the $d$-functions is to use the recurrence relations 
for $d$-functions. This will be discussed in a separate publication.

%
%
%
%
\setcounter{equation}{0}
\section{Concluding Remarks}

  We have constructed $h$-symplecton in this article and 
investigated some of its properties. It has been seen that 
many properties of $ sl(2) $ symplecton are inherited to $h$-symplecton. 
Unfortunately, $h$-dependence of $h$-symplecton is absorbed in $ \sigma $, 
namely, twist element $ \F$, so that we can not see specific 
hypergeometric function for $h$-deformation. It will become clear 
what kind of hypergeometric functions are specific to $h$-deformed 
quantities if we obtain explicit form  of $d$-function for 
$ SL_h(2) $ as in the case of $q$-deformed $ SU(2) $\cite{rep}.  
The $ sl(2) $ symplecton has a simple generating function.  
We presented (\ref{hgenfun}) as a generating function for 
$h$-symplecton. However, this may be one of possible choices, we 
might find simpler generating function. The use of quantum 
$h$-plane $ \xi_h, \eta_h $ instead of $ \xi, \eta $ is one of 
the possibilities. We have done some calculation to find 
simpler form of generating function in terms of 
$ \xi_h $ and $ \eta_h $, however, all what we 
obtained have more complicated form. 

  We would like to emphasize the usefulness of Lemma \ref{lem3}. 
This provides us a much simpler procedure to obtain $h$-symplecton 
than starting with the definition (\ref{defhsym}) and using the 
lemma, we could easily find another irrep. bases (\ref{hphix}) for 
$ SL_h(2). $ This lemma is, of course, applicable to any 
Jordanian quantum algebra, since we usually know the explicit form 
of twist element. Furthermore, the lemma is extended to quasitriangular 
Hopf algebras \cite{fiore2}. For quasitriangular Hopf algebras, the 
twist elements are usually not known, they are known up to certain order 
of the deformation parameters. It is expected that many properties 
of tensor operators for quasitriangular Hopf algebras are studied 
based on the present knowledge of the tensor operators for 
Lie algebras via Lemma \ref{lem3}, even if the 
explicit form of tensor operators is not obtained. It may also be possible to 
apply Lemma \ref{lem3} to the investigation of $q$-symplecton.

\section*{Appendices}
\appendix

%
%
%
%
\section{Relation to Ohn's $\sl2$}

Ohn defined in Ref.\cite{ohn} $ \sl2 $ as an algebra generated 
by $ H, X $ and $ Y $ subject to
\bea
  & & [X,\; Y] = H, \qquad 
      [H,\; X] = 2 {\displaystyle \frac{\sinh hX}{h}}, \label{ohnsl2}
  \nn \\
  & & [H,\; Y] = -Y (\cosh hX) - (\cosh hX) Y.
\eea
Meanwhile, the commutation relations of $ J_{\pm}, J_0 $, which are generators 
of $ \sl2 $ in this article, are same as $ sl(2) $. 
These two kinds of generators are related by 
\bea
  & & H = e^{-\sigma/2} J_0, 
      \qquad
      X = {\displaystyle \frac{\sigma}{2h}}, \label{relation} \\
  & & Y = e^{-\sigma/2} (J_- + {\displaystyle \frac{h}{2}J_0^2} ) 
      - {\displaystyle \frac{h}{8}} e^{\sigma/2} (e^{-\sigma}-1). \nn
\eea
By this relation, not only the commutation relations but also the Hopf algebra 
mappings are transformed each other. 
The relation (\ref{relation}) corresponds to the one parameter case discussed in Ref.\cite{na3} 
where two parameter Jordanian deformation of $ gl(2) $ is considered.

%
%
%
%
\setcounter{equation}{0}
\section{Hopf Algebra Structure of $ \sl2 $}

  We here give explicit formulae for the coproduct, counit 
and antipode of $ \sl2 $ calculated from (\ref{defHopf}).

\noindent
(i) coproduct
\bea
 \tilde \Delta(J_0) &=& J_0 \otimes e^{\sigma} + 1 \otimes J_0 \nn \\
 \tilde \Delta(J_+) &=& J_+ \otimes 1 + e^{-\sigma} \otimes J_+ \label{copsl2} \\
 \tilde \Delta(J_-) &=& J_- \otimes e^{\sigma} + 1 \otimes J_- 
    - h J_0 \otimes e^{\sigma} J_0 
    - {\textstyle \frac{h}{2}} J_0 (J_0+2) \otimes e^{\sigma} ( e^{\sigma}-1 ). \nn
\eea

\noindent
(ii) counit
\beq
  \tilde \epsilon(X) = 0, \qquad X = J_{\pm}, J_0. \label{cousl2}
\eeq 

\noindent
(iii) antipode
\bea
  & & \tilde S(J_0) = -J_0 e^{-\sigma}, \nn \\
  & & \tilde S(J_+) = - J_+ e^{\sigma}, \label{antisl2} \\
  & & \tilde S(J_-) = - J_- e^{-\sigma} -  {\textstyle \frac{h}{2}} J_0^2 
      (e^{-\sigma}+1) e^{-\sigma} + h J_0 (e^{-\sigma} -1) e^{-\sigma}. \nn
\eea
All of these are reduced to the ones for $ sl(2) $ in the limit of $ h=0 $. 
The Hopf algebra mappings for $ \sigma $ have  simple form
\beq
  \tilde \Delta(\sigma) = \sigma \otimes 1 + 1 \otimes \sigma, 
  \qquad
  \tilde \epsilon(\sigma) = 0, 
  \qquad
  \tilde S(\sigma) = - \sigma.       \label{sigmaHopf}
\eeq

%
%
%
%
\setcounter{equation}{0}
\section{Matrix Elements of $\F$}

  In this appendix, we show the explicit formula of matrix elements of the twist 
element $ \F $ (\ref{sigma}) and some of their properties. 
We denote a irrep. basis of $ \sl2 $ by the bracket notation $ \ket{jm} $ 
for the sake of simplicity. 

  It is easily verified the following relations from (\ref{defHopf})
\bea
 \tilde \Delta(J_{\pm})  \F \ket{j_1 m_1} \otimes \ket{j_2 m_2}
   &=& \sqrt{(j_1\mp m_1) (j_1 \pm m_1 +1)} \F \ket{j_1 \; m_1\pm 1} \otimes \ket{j_2 m_2}
 \nn \\
 &+& \sqrt{(j_2 \mp m_2) (j_2 \pm m_2 +1)} \F \ket{j_1 m_1} \otimes \ket{j_2 \; m_2\pm 1},
 \label{toCGC} \\
 \tilde \Delta(J_0)  \F \ket{j_1 m_1} \otimes \ket{j_2 \ m_2}
 &=&
 2(m_1 + m_2) \F \ket{j_1 m_1} \otimes \ket{j_2 m_2}. \nn
\eea
It shows that the vectors $ \F \ket{j_1 m_1} \otimes \ket{j_2 m_2} $ 
for $ \sl2 $ play the 
same role as $ \ket{j_1 m_1} \otimes \ket{j_2 m_2} $ for $ sl(2) $. Eq.(\ref{tensor}) 
is readily obtained from this. 
Another proof of Lemma \ref{lem1} with the bases of $ \sl2 $ in Ref.\cite{ohn} is 
found in Refs.\cite{na2, vdj}.

 In the bracket notation, matrix elements of $ \F $ are defined by 
\[
 F^{\ \ j_1,j_2}_{k_1,k_2\ m_1,m_2} = \bra{j_1 k_1} \otimes \bra{j_2 k_2} 
  {\cal F} \ket{j_1 m_1} \otimes \ket{j_2 m_2}. 
\]
We first show a relationship between the matrix elements of $ \F $ and 
its inverse

\beq
  F^{\ \ j_1,j_2}_{-n_1,-n_2 \ -m_1,-m_2} = (F^{-1})^{\ \ j_1,j_2}_{m_1,m_2 \ n_1,n_2}.
  \label{inv}
\eeq
The LHS of (\ref{inv}) is calculated as
\bea
  {\rm LHS} &=& \bra{j_1 -n_1}\otimes \bra{j_2 -n_2} 
     \sum_{\ell = 0}^{\infty} \frac{(-J_0)^{\ell} \otimes \sigma^{\ell}}
     {2^{\ell} \ell !} \ket{j_1 -m_1} \otimes \ket{j_2 -m_2} \nn \\
      &=& \delta_{n_1,m_1} \bra{j_2 -n_2} \exp(n_1 \sigma) \ket{j_2 -m_2}, 
     \label{lhs}
\eea
where $ \bra{j_1 -n_1} J_0 = -2n_1 \bra{j_1 -n_1} $ is used.
While the RHS is
\beq
  {\rm RHS} =  
  \bra{j_1 m_1} \otimes \bra{j_2 m_2} \sum_{\ell = 0}^{\infty} 
  \frac{J_0^{\ell} \otimes \sigma^{\ell}}
  {2^{\ell} \ell !} \ket{j_1 n_1} \otimes \ket{j_2 n_2} 
  = \delta_{m_1,n_1} \bra{j_2 m_2} \exp(n_1 \sigma) \ket{j_2 n_2}. 
  \label{rhs}
\eeq
Note that
\bea
 & & J_+ \ket{jm}= \sqrt{(j-m)(j+m+1)} \ket{j m+1}, \nn \\
 & & \bra{j -m} J_+ = \sqrt{(j-m)(j+m+1)} \bra{j -m-1}.\nn
\eea
It follows that any polynomials in $ J_+ $, denoted by $ f(J_+) $, 
satisfies
\beq
  \bra{jm}f(J_+)\ket{jn} = \bra{j -n} f(J_+) \ket{j -m}.
  \label{poly}
\eeq
Since $ \sigma $ is a polynomial in $ J_+ $, we see that the 
(\ref{lhs}) equals to (\ref{rhs}). Thus (\ref{inv}) has been proved.

 We next show that the matrix elements of $ \F $ are given by
\bea
  & & F^{\quad j_1,j_2}_{k_1, k_2\ m_1, m_2} = 
  \delta_{k_1, m_1} \theta(m_2 \leq k_2 \leq j_2) \; S^{\ j_2}_{k_2,m_2}
  \nn \\
  & & \qquad \times 
  \left\{
  \begin{array}{ll}
   \displaystyle{
   \frac{(2k_2 - 2m_1 -2m_2-2)!!}{(k_2-m_2)! (-2m_1-2)!!} \; h^{k_2-m_2}_,
   }
   & {\rm for} \ \; m_1 \leq 0 \\
    \\
   (-2h)^{k_2-m_2}  
    \displaystyle{
     \sum_{\ell=0}^{j_2-m_2} (-1)^{\ell} \left(
       \begin{array}{c}
         2m_1 \\ k_2-m_2-\ell
       \end{array}
       \right)
     \frac{(2\ell + 2m_1-2)!!}{2^{\ell} \ell ! (2m_1-2)!!}_,
   }
   & {\rm for}\ \; m_1 > 0 
  \end{array}
  \right.
  \label{emF}
\eea
where  $ n !! = 1 $ for $ n \leq 0 $ 
and $ \theta(m_2 \leq k_2 \leq j_2) =1 $ if and only if  the inequality in the parenthesis 
holds, otherwise $ \theta $ vanishes.  $ S^{\ j_2}_{k_2,m_2} $ is defined by
\[
  S^{\ j_2}_{k_2,m_2} = \left\{
  \frac{(j_2-m_2)! (j_2 + k_2)!}{(j_2+m_2)! (j_2 - k_2)!}
  \right\}^{1/2}_.
\]
To prove (\ref{emF}), note that similar to (\ref{rhs}) we have
\beq
  F^{\quad j_1,j_2}_{k_1, k_2\ m_1, m_2} = \delta_{k_1,m_1} \bra{j_2 k_2} e^{-m_1 \sigma} 
  \ket{j_2 m_2}. 
  \label{emF2}
\eeq
One can use the power series expansion in order to compute the RHS of (\ref{emF2})
\beq
  (1-X)^{-\ell/2} = \sum_{n=0}^{\infty}\; 
  \frac{(2n+\ell-2)!!}{2^n n! (\ell-2)!!} X^n, 
  \qquad
  \ell \in {\bf Z}_{+}.                 \label{power}
\eeq

\noindent
(i) For $ m_1 \leq 0 $

  Let $ m_1 = -\ell/2 \ ( \ell \in {\bf Z}_{+}) $ and using (\ref{power})
\bea
 e^{-m_1 \sigma} \ket{j_2 m_2} &=& 
 (1-2hJ_+)^{-\ell/2} \ket{j_2 m_2} \nn \\  
  &=& 
  \sum_{n=0}^{j_2-m_2}\; \frac{(2n+\ell-2)!!}{2^n n! (\ell-2)!!} (2h)^n 
  \left\{
  \frac{(j_2-m_2)! (j_2 + m_2 +n)!}{(j_2+m_2)! (j_2 - m_2-n)!}
  \right\}^{1/2} \ket{j_2\; m_2+n}. \nn
\eea
Therefore $ \bra{j_2 k_2} e^{-m_1 \sigma} \ket{j_2 m_2} $ takes values if and 
only if $ k_2 = m_2 + n $. This proves the first part of (\ref{emF}). 

\noindent
(ii) For $ m_1 > 0 $

 Let $ m_1 = \ell/2 \ ( \ell \in {\bf Z}_{+}). $ Since
\[
  e^{-m_1 \sigma} \ket{j_2 m_2} = e^{-\ell \sigma} e^{\ell \sigma /2} \ket{j_2 m_2},
\]
we can apply the previous result to compute $ e^{\ell \sigma /2} \ket{j_2 m_2} $ and then 
applying the binomial expansion to $ e^{-\ell \sigma} = (1-2hJ_+)^{\ell} $
\bea
  & & e^{-m_1 \sigma} \ket{j_2 m_2} \nn \\
  & & = 
  \sum_{t=0}^{j_2-m_2}\; \frac{(2t+\ell-2)!!}{2^t t! (\ell-2)!!} (2h)^t 
  \left\{
  \frac{(j_2-m_2)! (j_2 + m_2 +t)!}{(j_2+m_2)! (j_2 - m_2-t)!}
  \right\}^{1/2} e^{-\ell \sigma} \ket{j_2\; m_2+t}   
  \nn \\
  & & = 
  \sum_{n=0}^{\ell} \; \sum_{t=0}^{j_2-m_2}\; 
  \left(
       \begin{array}{c}
         \ell \\ n
       \end{array}
  \right)
  \frac{(2t+\ell-2)!!}{2^t t! (\ell-2)!!} (-1)^n (2h)^{t+n} 
  \left\{
  \frac{(j_2-m_2)! (j_2 + m_2 +t+n)!}{(j_2+m_2)! (j_2 - m_2-t-n)!}
  \right\}^{1/2}
  \nn \\
  & & \ \ \times \ket{j_2\; m_2+t+n}. \nn
\eea
Replacing $ t+n $ with $n$, we obtain
\bea
  e^{-m_1 \sigma} \ket{j_2 m_2} 
  &=&
  \sum_{t,n}\;  (-1)^t (-2h)^n \left(
       \begin{array}{c}
         \ell \\ n-t
       \end{array}
  \right) 
  \nn \\
  &\times & 
  \frac{(2t+\ell-2)!!}{2^t t! (\ell-2)!!} 
  \left\{
  \frac{(j_2-m_2)! (j_2 + m_2 +n)!}{(j_2+m_2)! (j_2 - m_2-n)!}
  \right\}^{1/2} \ket{j_2\; m_2+n}. \nn
\eea
Again  $ \bra{j_2 k_2} e^{-m_1 \sigma} \ket{j_2 m_2} $ takes values if and 
only if $ k_2 = m_2 + n $. This completes the proof of (\ref{emF}). 

 We can obtain the explicit formula for the universal $R$-matrix in the irreps. with 
highest weight $ j_1 $ and $ j_2 $ by combining the (\ref{emF}) and 
relation (\ref{inv}), since universal $R$-matrix for $ \sl2 $ is given by 
$ {\cal R} = \F_{21} \F^{-1}. $

%
%
%
%
\setcounter{equation}{0}
\section{Proof of Lemma \ref{lem2}}

 Let $ V^a, V^b $ and $ V^c $ be representation spaces of $ \sl2 $ 
with highest weight $ a, b $ and $c$, respectively. Bases of each 
space are denoted as $ e^a_{\alpha},\ -a \leq \alpha \leq a. $ 
We would like to construct irrep. bases in the space 
$ V^a \otimes V^b \otimes V^c $ in two ways, namely, 
$ (V^a \otimes V^b) \otimes V^c $ and $ V^a \otimes (V^b \otimes V^c) $. 
According to the discussion in Appendix B, irrep. bases in the space 
$ V^a \otimes V^b $ are given by
\[
  e^{(ab)d}_{\delta} = \sum C^{a,b,d}_{\alpha, \beta, \delta} \; \F \; 
  e^a_{\alpha} \otimes e^b_{\beta}.
\]
Then we couple these with the bases in $ V^c $ to obtain
\[
  \psi^e_{\epsilon} = \sum C^{d,c,e}_{\delta, \gamma, \epsilon} \; 
  (\tilde \Delta \otimes id)(\F) \; e^{(ab)d}_{\delta} \otimes c^c_{\gamma} 
  = 
  \sum C^{d,c,e}_{\delta, \gamma, \epsilon} \; C^{a,b,d}_{\alpha, \beta, \delta} \; 
   (\tilde \Delta \otimes id)(\F) \F_{12} \; 
  e^a_{\alpha} \otimes e^b_{\beta} \otimes e^c_{\gamma}.
\]
Similarly we obtain the following bases when we couple $ V^b $ and $ V^c $ first
\[
  \psi'^e_{\epsilon} = \sum C^{a,f,e}_{\alpha, \rho, \epsilon} \; C^{b,c,f}_{\beta, \gamma, \rho} \; 
  (id \otimes \tilde \Delta)(\F) \F_{23} \;
  e^a_{\alpha} \otimes e^b_{\beta} \otimes e^c_{\gamma}.
\]
From (\ref{sigmaHopf})
\[
  (id \otimes \tilde \Delta)(\F) = \exp(-{\textstyle \frac{1}{2}} 
  J_0 \otimes \tilde \Delta(\sigma)) = \F_{12} \F_{13}.
\]
Using the relations (\ref{defHopf}), (\ref{twist2}) and above
\bea
  (\tilde \Delta \otimes id)(\F) 
   &=& \F_{12}(\Delta \otimes id) (\F) \F^{-1}_{12} 
    = \F_{23} (id \otimes \Delta) (\F) \F^{-1}_{12} \nn \\
   &=& (id \otimes \tilde \Delta) (\F) \F_{23} \F^{-1}_{12} \nn \\
   &=& \F_{12} \F_{13} \F_{23} \F^{-1}_{12}. \nn
\eea
It follows that $  \psi^e_{\epsilon} $ and $  \psi'^e_{\epsilon} $ are 
rewritten as
\bea
 & &  \psi^e_{\epsilon} = 
 \sum C^{d,c,e}_{\delta, \gamma, \epsilon} \; C^{a,b,d}_{\alpha, \beta, \delta} \; 
 \F_{12} \F_{13} \F_{23} \; e^a_{\alpha} \otimes e^b_{\beta} \otimes e^c_{\gamma},
 \nn \\
 & & \psi'^e_{\epsilon} = 
 \sum C^{a,f,e}_{\alpha, \rho, \epsilon} \; C^{b,c,f}_{\beta, \gamma, \rho} \; 
 \F_{12} \F_{13} \F_{23} \; e^a_{\alpha} \otimes e^b_{\beta} \otimes e^c_{\gamma}. 
 \nn
\eea
The Racha coefficients $ W_h(abce ; df) $ for $ \sl2 $ is defined by
\beq
  \psi^e_{\epsilon} = \sum_f \; \sqrt{(2d+1) (2f+1)} W_h(abce ; df) \; \psi'^e_{\epsilon}.
  \label{hRacha}
\eeq
It it now obvious that the Racha coefficients for $ \sl2 $ satisfy the relation
\beq
 \sum_{\delta} \; C^{d,c,e}_{\delta, \gamma, \epsilon} \; C^{a,b,d}_{\alpha, \beta, \delta} \; 
 = 
 \sum_{f, \rho}\; C^{a,f,e}_{\alpha, \rho, \epsilon} \; C^{b,c,f}_{\beta, \gamma, \rho} \; 
 \sqrt{(2d+1) (2f+1)} W_h(abce ; df). 
 \label{recRacha}
\eeq
This is the same relation for the Racha coefficients for $ sl(2) $. 
This proves Lemma  \ref{lem2}.

%
%
%
%


\begin{thebibliography}{99}
\bibitem{BieL8} See for example, L. C. Biedenharn and J. D. Louck, 
\textit{Angular Momentum in Quantum Physics: Theory and Application, 
Encyclopedia of Mathematics and Its Applications} 
{\bf 8}, Addison-Wesley, Reading, Massachusetts (1981). 
\bibitem{BLohe1} See for example, L. C. Biedenharn and M. A. Lohe, 
\textit{Quantum Group Symmetry and $q$-Tensor Algebras}, 
World Scientific (1995).
\bibitem{dem} E. E. Demidov \textit{et al}, Prog. Theor. Phys. Suppl. 
{\bf 102},  203 (1990).
\bibitem{zak} S. Zakrewski, Lett. Math. Phys. {\bf 22},  287 (1991).
\bibitem{ewen} H. Ewen, O. Ogievetsky and J. Wess, Lett. Math. Phys. 
{\bf 22}, 297 (1991).
\bibitem{ohn} Ch. Ohn, Lett. Math. Phys. {\bf 25},  85 (1992).
\bibitem{BieL} L. C. Biedenharn and J. D. Louck, Ann. Phys. (N.Y.) {\bf 63}, 
 459 (1971).
\bibitem{BieL2} L. C. Biedenharn and J. D. Louck, \textit{The Racha-Wigner 
Algebra in Quantum Theory, Encyclopedia of Mathematics and Its Applications} 
{\bf 9}, Addison-Wesley, Reading, Massachusetts (1981).
\bibitem{LBL} M. A. Lohe, L. C. Biedenharn and J. D. Louck, Phys. Rev. {\bf D43}, 
 617 (1991).
\bibitem{BLbook} L. C. Biedenharn and M. A. Lohe, in \textit{Quantum Groups}, 
Proceedings of the Argonne Workshop, ed. T. Curtright, D. Fairlie and 
C. Zachos, World Scientific (1991).
\bibitem{Nom} M. Nomura, J. Math. Phys. {\bf 33},  3636 (1992).
\bibitem{GelF} I. M. Gelfand and D. B. Fairlie, Comm. Math. Phys. {\bf 136}, 
 487 (1991).
\bibitem{rs} V. Rittenberg and M. Scheunert, J. Math. Phys. {\bf 33}, 
436  (1992).
\bibitem {dri} V. G. Drinfeld, Leninglrad Math. J. {\bf 1}, 1419  (1990). \\
See also, N. Yu. Reshetikhin, Lett. Math. Phys. {\bf 20}, 331  (1990).
\bibitem{ks} P. P. Kulish and A. A. Stolin, Czech. J. Phys. {\bf 47}, 
 123 (1997).
\bibitem{frt} L. D. Faddeev, N. Yu. Reshetikhin and L. A. Takhtajan, 
Leningrad Math. J. {\bf 1},  193 (1990).
\bibitem{na1} N. Aizawa, J. Phys. {\bf A}:Math. Gen. {\bf 31},  5467 (1998).
\bibitem{ques1} C. Quesne, Czech. J. Phys. {\bf 48}, 1471  (1998).
\bibitem{fiore} G. Fiore, J. Math. Phys. {\bf 39},  3437 (1998) .
\bibitem{kar} V. Karimipour, Lett. Math. Phys. {\bf 30}, 87  (1994).
\bibitem{CQ} R. Chakrabarti and C. Quesne, preprint, math.QA/9811064.
\bibitem{rep} T. Masuda $ et. al. $, J. Funct. Anal. {\bf 99},  357 (1991), and 
references therein.
\bibitem{fiore2} G. Fiore, preprint, q-alg/9708017.
\bibitem{na3} N. Aizawa, Czech. J. Phys. {\bf 48},  1273 (1998).
\bibitem{na2} N. Aizawa, J. Phys. {\bf A}:Math. Gen. {\bf 30},  5981 (1997).
\bibitem{vdj} J. Van der Jeugt, J. Phys. {\bf A}:Math. Gen {\bf 31},  1495 (1998), \ \ 
Czech. J. Phys. {\bf 47},  1283 (1997).
\bibitem{btoq} P. P. Kulish and E. V. Damaskinsky, J. Phys. {\bf A}:Math. Gen. 
{\bf 23}, L415  (1990). \\
M. Chaichian, P. P. Kulish and J. Lukierski, Phys. Lett. {\bf B262},  43 (1991).
\end{thebibliography}
\end{document}